\documentclass{birkjour}
 \newtheorem{thm}{Theorem}[section]
 \newtheorem{cor}[thm]{Corollary}
 \newtheorem{lem}[thm]{Lemma}
 
 \theoremstyle{definition}
 
 \theoremstyle{remark}
 \newtheorem{rem}[thm]{Remark}
 
 \numberwithin{equation}{section}

\newcommand{\R}{\mathbb{R}}

\newcommand{\N}{\mathbb{N}}

\newcommand{\T}{{\mathrm{T}}}

\newcommand{\cont}{{\mathrm C }}

\newcommand{\contX}{\cont(X)}

\newcommand{\orderdual}[1]{{#1}^\sim}

\newcommand{\ordercontc}[1]{{#1}^\sim_{\mathrm c}}
\newcommand{\zerof}{{\mathbf 0}}
\newcommand{\onef}{{\mathbf 1}}

\newcommand{\spn}{\mathrm{ span}}

\newcommand{\Zsets}[1]{\mathbf{Z}({#1})}
\newcommand{\ZX}{\Zsets{X}}
\newcommand{\CZsets}[1]{\mathbf{Z}^c({#1})}
\newcommand{\CZX}{\CZsets{X}}

\newcommand{\nbh}[1]{\mathcal{N}_{{#1}}}
\newcommand{\cnbh}[1]{\mathcal{N}_{{#1}}^\ast}
\newcommand{\znbh}[1]{\mathcal{N}^z_{{#1}}}
\newcommand{\cznbh}[1]{\mathcal{N}^c_{{#1}}}

\begin{document}

%-------------------------------------------------------------------------
% editorial commands: to be inserted by the editorial office
%
%\firstpage{1} \volume{228} \Copyrightyear{2004} \DOI{003-0001}
%\seriesextraline{This is the Concrete Title of this Book\br H.E. R and S.T.C. W, Eds.}
%\firstpage{1}

\title[Universally complete $\contX$]
 {Universally complete spaces of continuous functions}

%----------Author 1
\author[J. H. van der Walt]{Jan Harm van der Walt}

\address{%
Department of Mathematics and Applied Mathematics\\ 
University of Pretoria\\
Cor\-ner of Lynnwood Road and Roper Street\\
Hatfield 0083, Pretoria\\
   South Africa}

\email{janharm.vanderwalt@up.ac.za}



\date{August 24, 2020}

%%% ----------------------------------------------------------------------

\begin{abstract}
We characterise Tychonoff spaces $X$ so that $\contX$ is universally $\sigma$-complete and universally complete, respectively.
\end{abstract}

%%% ----------------------------------------------------------------------
\maketitle
%%% ----------------------------------------------------------------------
%\tableofcontents

\section{Introduction}\label{Sec:  Introduction}

Recently, Mozo Carollo \cite{MozoCarollo2020} showed, in the context of point-free topology, that the vector lattice $\contX$ of continuous, real valued functions on a Tychonoff (completely regular $\T_1$) space $X$ is universally complete if and only if $X$ is an extremally disconnected P-space.  This paper aims to make this result and its proof accessible to those members of the positivity community who, like the author, are less familiar with point-free topology.  In so doing, and based on results due to Fremlin \cite{Fremlin1975} and Veksler and Ge\u{\i}ler \cite{VekslerGeiler1972}, we obtain a refinement of Mozo Carollo's result.  In particular, we characterise those Tychonoff spaces $X$ for which $\contX$ is laterally $\sigma$-complete.  We also include some remarks on $\sigma$-order continuous duals of spaces $\contX$ which are universally $\sigma$-complete.

The paper is organised as follows.  In Section \ref{Sec:  Preliminaries} we introduce definitions and notation used throughout the paper, and recall some results from the literature.  Section \ref{Sec:  Universally complete C(X)} contains the main results of the paper, namely, characterisations of those Tychonoff spaces $X$ for which $\contX$ is universally complete and universally $\sigma$-complete, respectively.  %The final section contains some examples.

\section{Preliminaries}\label{Sec:  Preliminaries}

Throughout this paper $X$ denotes a Tychonoff space; that is, a completely regular $\T_1$ space.  $\contX$ stands for the lattice of all real-valued and continuous functions on $X$.  For $u\in\contX$, $Z(u)$ denotes the zero set of $u$; that is, $Z(u)=u^{-1}[\{0\}]$.  The co-zero set of $u$ is $Z^c(u) = X\setminus Z(u)$.  The collection of zero sets in $X$ is denoted $\ZX$, while $\CZX$ consists of all co-zero sets in $X$.  For $x\in X$ the collection of open neighbourhoods of $x$ is denoted $\nbh{x}$, and $\cnbh{x}$ denotes the set of clopen neighbhourhoods of $x$.  A zero-neighbourhood of $x\in X$ is a set $V\in \ZX$ so that $x$ belongs to the interior of $V$.  The collection of all zero-neighbourhoods of $x\in X$ is denoted $\znbh{x}$, and $\cznbh{x}=\nbh{x}\cap \CZX$.  Observe that $\CZX$ is a basis for the topology of $X$.  Hence for each $x\in X$ and every $V\in \nbh{x}$ there exists $U\in \cznbh{x}$ so that $U\subseteq V$.  Furthermore, for every $V\in \nbh{x}$ there exists $W\in \znbh{x}$ so that $W\subseteq V$.  The standard reference for all of this is \cite{GillmanJerison1960}.

We write $\onef$ for the function which is constant one on $X$.  More generally, for $A\subseteq X$, the indicator function of $A$ is $\onef_A$.  The constant zero function is $\zerof$.

We recall, see for instance \cite{GillmanJerison1960}, that $X$ is \begin{itemize}
%    \item[(i)] \emph{zero-dimensional}\footnote{The term zero-dimensional should be understood in terms of small inductive dimension \cite[Definition 1.1.1 \& Proposition 1.2.1]{Engelking}, as opposed to the Lebesgue covering dimension used in \cite{GillmanJerison1960}.} if it has a basis consisting of clopen sets;
    \item[(i)] \emph{basically disconnected} if the closure of every co-zero set is open;
    \item[(ii)]\emph{extremally disconnected} if the closure of every open set is open.
\end{itemize}
Every extremally disconnected space is basically disconnected, but not conversely \cite[Problem 4N]{GillmanJerison1960}.  Since $\CZX$ is a basis for the topology on $X$, every basically disconnected space is zero-dimensional\footnote{The term zero-dimensional should be understood in terms of small inductive dimension \cite[Definition 1.1.1 \& Proposition 1.2.1]{Engelking}, as opposed to the Lebesgue covering dimension used in \cite{GillmanJerison1960}.}; that is, it has a basis consisting of clopen sets.  The converse is false.  For instance, $\mathbb{Q}$ is zero-dimensional, the set of all open intervals with irrational endpoints forming a basis of clopen sets, but not basically disconnected, since $(0,1)$ is a co-zero set whose closure is not open.

Each of the properties (i) and (ii) of $X$ corresponds to order-theoretic properties of $\contX$, see for instance \cite[Theorems 43.2, 43.3, 43.8 \& 43.11]{LuxemburgZaanen1971RSI}.  In particular, $X$ is \begin{itemize}
%    \item[(i*)] \emph{zero-dimensional} if and only if ...$\contX$ has sufficiently many projections;
    \item[(i*)] basically disconnected if and only if $\contX$ is Dedekind $\sigma$-complete, if and only if $\contX$ has the principle projection property;
    \item[(ii*)] extremally disconnected if and only if $\contX$ is Dedekind complete, if and only if $\contX$ has the projection property.
\end{itemize}

$X$ is a P-space \cite{GillmanHenriksenTAMS1954} if the intersection of countably many open sets in $X$ is open.  Equivalently, $X$ is a P-space if $Z(u)$ is open (hence clopen) for every $u\in\contX$.  Clearly, every discrete space is a P-space, but the converse is false, see \cite[Problem 4N]{GillmanJerison1960}.  In fact, there exists a P-space without any isolated points \cite[Problem 13P]{GillmanJerison1960}.  Evidently, every P-space is basically disconnected (in particular, every $Z\in\ZX$ is open), but not conversely, see \cite[Problem 4M]{GillmanJerison1960}.

The following basic lemma may well be known, but we have not found it in the literature.  We include the simple proof for the sake of completeness.

\begin{lem}\label{Lemma:  L1}
Let $X$ be zero-dimensional.  Then the following statements are equivalent. \begin{itemize}
    \item[(i)] $X$ is a P-space.
    \item[(ii)] The intersection of countably many clopen sets is clopen.
    \item[(iii)] The union of countably many clopen sets is clopen.
\end{itemize}
\end{lem}
\begin{proof}
By definition, (i) implies (ii) and (iii), and, (ii) and (iii) are equivalent.  It therefore suffices to show that (ii) implies (i).

Assume that (ii) is true.  For each $n\in \mathbb{N}$ let $U_n$ be an open subset of $X$.  Let $U=\displaystyle \bigcap\{U_n ~:~ n\in\N\}$.  If $U=\emptyset$ we are done, so assume that $U\neq \emptyset$. Fix any $x\in U$.  Since $X$ is zero-dimensional, there exists for each $n\in\mathbb{N}$ a set $V_n\in\cnbh{x}$ so that $V_n\subseteq U_n$.  Let $V=\bigcap_{n\in\mathbb{N}} V_n$.  Then $x\in V\subseteq U$ and, by assumption, $V$ is clopen, hence open.  Therefore $U$ is open so that $X$ is a P-space.
\end{proof}

We recall, for later use, the following results of Fremlin \cite{Fremlin1975} and Veksler and Ge\u{\i}ler \cite{VekslerGeiler1972}, respectively; see also \cite{AliprantisBurkinshaw78}.

\begin{thm}\label{Theorem:  T1}
Let $L$ be a Dedekind complete vector lattice.  Then the following statements are equivalent. \begin{itemize}
    \item[(i)] $L$ is universally complete.
    \item[(ii)] $L$ is universally $\sigma$-complete and has a weak order unit.
\end{itemize}
\end{thm}

%\begin{thm}\label{Theorem:  T7}
%Let $L$ be a universally $\sigma$-complete vector lattice. Then every positive linear functional $\varphi:L\to\R$ is a finite linear combination of real-valued Riesz homomorphisms on $L$.
%\end{thm}

\begin{thm}\label{Theorem:  T2}
Let $L$ be an Archimedean vector lattice.  The following statements are true. \begin{itemize}
    \item[(i)] If $L$ is laterally complete then $L$ has the projection property.
    \item[(ii)] If $L$ is laterally $\sigma$-complete then $L$ has the principle projection property.
\end{itemize}
\end{thm}

\section{Universally complete $\contX$}\label{Sec:  Universally complete C(X)}

We begin this section with a characterisation of those $X$ for which $\contX$ is universally $\sigma$-complete.

\begin{thm}\label{Theorem:  T5}
The following statements are equivalent. \begin{itemize}
    \item[(i)] $\contX$ is laterally $\sigma$-complete.
    \item[(ii)] $\contX$ is universally $\sigma$-complete.
    \item[(iii)] $X$ is a P-space.
\end{itemize}
\end{thm}

\begin{proof}
Assume that $\contX$ is laterally $\sigma$-complete.  It follows from Theorem \ref{Theorem:  T2} (ii) that $\contX$ has the principle projection property.  Therefore $\contX$ is Dedekind $\sigma$-complete, hence universally $\sigma$-complete.  Conversely, if $\contX$ is universally $\sigma$-complete then, by definition, it is laterally $\sigma$-complete.  Hence (i) and (ii) are equivalent.

Assume that $\contX$ is laterally $\sigma$-complete.  Then $\contX$ has the principle projection property so that $X$ is basically disconnected, hence zero-dimensional.  We show that $X$ is a P-space.

By Lemma \ref{Lemma:  L1} it suffices to show that the intersection of countably many clopen subsets of $X$ is clopen. Assume that $U_k\subseteq X$ is clopen for each $k\in \mathbb{N}$, and let $U=\bigcap\{ U_k ~:~ k\in\mathbb{N}\}$.  We claim that $U$ is clopen.

Let $V_0=X$, $V_1=U_1$ and, for each natural number $n>1$, let $V_n = U_1 \cap\ldots\cap U_n$.  Then $V_n$ is clopen for each $n\in\mathbb{N}$, $U= \bigcap\{V_n ~:~ n\in\mathbb{N}\}$ and $V_{n+1}\subseteq V_n$ for every $n\in\mathbb{N}$.  For each $n\in\mathbb{N}$, let $W_n =V_{n-1}\setminus V_n$.  Then each $W_n$ is clopen and $W_n\cap W_m=\emptyset$ whenever $n\neq m$. Moreover, $\bigcup\{W_n ~:~ n\in\mathbb{N}\}=X\setminus U$.  Indeed, the inclusion $\bigcup\{ W_n ~:~ n\in\mathbb{N}\} \subseteq X\setminus U$ is immediate. For the reverse inclusion, consider some $x\in X\setminus U$.  There exists $n\in\mathbb{N}$ so that $x\in X\setminus V_n$.  Let $n_0 = \min\{ n\in\mathbb{N} ~:~ x\in X\setminus V_n\}$.  Then, since $V_0=X$, $x\in V_{n-1}\setminus V_n=W_n$.  Hence $x\in \bigcup\{W_n ~:~ n\in\mathbb{N}\}$.

Let $w_n = n \onef_{W_n}$, $n\in\mathbb{N}$, and $F=\{w_n ~:~ n\in\mathbb{N}\}$. Then $F\subseteq \contX^+$ and the $w_n$ are mutually disjoint.  Therefore, since $\contX$ is universally $\sigma$-complete, $w=\sup F$ exists in $\contX$.

Fix $x\in U$.  There exists $V\in\nbh{x}$ so that $w(y)<w(x)+1$ for all $y\in V$.  Fix a natural number $N_0\geq w(x)+1$. Then, for all $n\geq N_0$ and $y\in W_n$, $w(y)\geq w_n (y) = n \geq N_0 \geq w(x)+1$ so that $y\notin V$.  Therefore $V\cap W_n = \emptyset$ for all $n\geq N_0$.  Let $W= V_{N_0}\cap V$.  Then $W\in\nbh{x}$ and, since $W_n\cap V_{N_0}=\emptyset$ for all $n<N_0$, $W\cap W_n=\emptyset$ for all $n\in\mathbb{N}$.  Therefore $W\subseteq X\setminus \bigcup\{ W_n ~:~ n\in\mathbb{N}\}=U$.  This shows that $U$ is open, and, since each $U_k$ is closed, $U$ is also closed, hence clopen.  By Lemma \ref{Lemma:  L1}, $X$ is a P-space.  Hence (i) implies (iii)

Assume that $X$ is a P-space.  Consider a countable set $F$ of mutually disjoint elements of $\contX^+$.  We observe that for each $x\in X$ there is at most one $u\in F$ so that $u(x)>0$.  Hence the function
\[
w:X\ni x\to \sup\{u(x) ~:~ u\in F\}\in\mathbb{R}^+
\]
is well defined.  We claim that $w\in \contX$ so that $w=\sup F$ in $\contX$.

Fix $x\in X$.  Assume that $w(x)>0$.  Then there exists $u\in F$ and $V\in \nbh{x}$ so that $u(y)=w(y)>0$ for all $y\in V$.  Hence $w$ is continuous at $x$.  Suppose $w(x)=0$.  Then $u(x)=0$ for all $u\in F$.  Since $X$ is a P-space there exists for each $u\in F$ some $V_u\in \nbh{x}$ so that $u(y)=0$ for every $y\in V_u$.  The set $V=\cap\{V_u ~:~ u\in F\}$ is an open neighbourhood of $x$, and $w(y)=0$ for all $y\in V$.  Hence $w$ is continuous at $x$.  Thus $w$ is continuous at every $x\in X$, hence on $X$.  Therefore $\contX$ is laterally $\sigma$-complete.  Hence (iii) implies (i).
\end{proof}

Mozo Carollo's characterisation of those $X$ for which $\contX$ is universally complete now follows easily.

\begin{cor}\label{Corollary:  C1}
The following statements are equivalent. \begin{itemize}
    \item[(i)] $\contX$ is laterally complete.
    \item[(ii)] $\contX$ is universally complete.
    \item[(iii)] $X$ is an extremally disconnected P-space.
\end{itemize}
\end{cor}

\begin{proof}
Assume that $\contX$ is laterally complete.  By Theorem \ref{Theorem:  T2} (i), $\contX$ has the projection property and is therefore Dedekind complete, hence universally complete.  Conversely, if $\contX$ is universally complete, then it is laterally complete.  Therefore (i) and (ii) are equivalent.

Assume that $\contX$ is universally complete.  Then, since $\contX$ is Dedekind complete, $X$ is extremally disconnected, and by Theorem \ref{Theorem:  T5}, $X$ is a P-space.

Suppose that $\contX$ is an extremally disconnected P-space.  Then $\contX$ is Dedekind complete and, by Theorem \ref{Theorem:  T5}, laterally $\sigma$-complete.  Since $\contX$ has a weak order unit, it is universally complete by Theorem \ref{Theorem:  T1}.
\end{proof}

\begin{rem}
Isbell \cite{Isbell1955} showed that if $X$ is an extremally disconnected P-space, and $X$ has non-measurable cardinal, then $X$ is discrete.  It is consistent with ZFC that every cardinal is non-measurable.  %On the other hand, if the existence of a measurable cardinal is consistent with ZFC, then so is the non-existence of such cardinals.  Furthermore, ZFC cannot establish the consistency of the existence of a measurable cardinal with ZFC, assuming the consistency of ZFC itself.
\end{rem}

Combining Theorem \ref{Theorem:  T5} with \cite[Theorem 10.2]{dePagterHuijsmans1980II} we have the following.

\begin{cor}
The following statements are equivalent. \begin{itemize}
    \item[(i)] $\contX$ is laterally $\sigma$-complete.
    \item[(ii)] $\contX$ is universally $\sigma$-complete.
    \item[(iii)] $X$ is a P-space.
    \item[(iv)] $\contX$ is a von Neumann regular\footnote{For every $u\in \contX$ there exists $v\in\contX$ so that $u=vu^2$.} ring.
    \item[(v)] $\contX$ is $z$-regular\footnote{Every proper prime $z$-ideal in $\contX$ is a minimal prime $z$-ideal, see \cite{dePagterHuijsmans1980I,dePagterHuijsmans1980II} for details.}.
\end{itemize}
\end{cor}

\begin{rem}\label{Remark: R1}
Recall that a space $X$ is called realcompact\footnote{Realcompact spaces were introduced by Hewitt \cite{Hewitt1948} under the name ``Q-spaces", and defined as follows:  $X$ is a Q-space if every free maximal ring ideal in $\contX$ is hyper-real.  See for instance \cite[Problem 8A no. 1]{GillmanJerison1960} for the equivalence of our definition and Hewitt's.} if for every Tychonoff space $Y$ containing $X$ as a proper dense subspace, the map $\cont(Y)\ni f\mapsto f|_X\in \contX$ is not onto; that is, $X$ is not C-embedded in $Y$, see \cite[page 214]{Engelking1989}

If $X$ is a realcompact P-space, then $\orderdual{\contX}$ has a peculiar structure.  Indeed, due to a result of Fremlin \cite[Proposition 1.15]{Fremlin1975}, every $\varphi \in \orderdual{\contX}$ is a finite linear combination of linear lattice homomorphisms form $\contX$ into $\R$.  Xiong \cite{Xiong1989} showed that every such homomorphism is a positive scalar multiple of a point evaluation.  Hence
\[
\orderdual{\contX} = \spn\{\delta_x ~:~ x\in X\}=c_{00}(X).
\]
However, each $\delta_x$ is $\sigma$-order continuous. Indeed, consider a decreasing sequence $(u_n)$ in $\contX^+$ so that $\inf\{u_n(x) ~:~ n\in\N\}>0$ for some $x\in X$.  Then there exists a real number $\epsilon>0$ so that for every $n\in\N$ there exists $V_n\in\nbh{x}$ such that $u_n(y)>\epsilon$ for every $y\in V_n$.  Since $X$ is a P-space, $V=\bigcap\{V_n ~:~ n\in\N\}$ is open.  Therefore there exists $v\in\contX$ so that $\zerof< v \leq \epsilon \onef$ and $v(y)=0$ for $y\in X\setminus V$.  Since $u_n(y)>\epsilon$ for all $y\in V$ and $n\in\N$ it follows that $\zerof \leq v\leq u_n$ for all $n\in\N$; hence $u_n$ does not decrease to $\zerof$ in $\contX$.  This shows that $\delta_x\in\ordercontc{\contX}$.

Combining all of the above, we see that
\[
\orderdual{\contX} = \spn\{\delta_x ~:~ x\in X\}=c_{00}(X) = \ordercontc{\contX}.
\]
\end{rem}

\begin{rem}
The condition that $\delta_x\in \ordercontc{\contX}$ for all $x\in X$ does not imply that $X$ is a P-space.  In fact, this property characterises the so called almost-P-spaces introduced by Veksler \cite{Veksler1973}, see also \cite{Levy1977}.  A space $X$ is an almost-P-space if the nonempty intersection of countably many open sets has nonempty interior; equivalently, every $Z\in\ZX$ has nonempty interior.  Thus every P-space is an almost-P-space, but not conversely, see \cite{Levy1977}.

De Pagter and Huijsmans \cite{dePagterHuijsmans1980II} showed that $\contX$ has the $\sigma$-order continuity property if and only if $X$ is an almost-P-space.  Hence, if $X$ is an almost-P-space, then $\delta_x\in\ordercontc{\contX}$ for every $x\in X$.  For the converse, suppose that $X$ is not an almost-P-space.  Then there exists $u\in\contX^+$ so that $Z(u)$ has empty interior.  For each $n\in\N$, let $u_n = \onef \wedge (nu)$.  Then $(u_n)$ is increasing and bounded above by $\onef$.  Let $v\in \contX$ be an upper bound for $(u_n)$.  If $x\in X\setminus Z(u)$ then $\sup \{u_n(x) ~:~n\in\mathbb{N}\}=1$ so that $v(x)\geq 1$.  Since $Z(u)$ has empty interior and $v$ is continuous, it follows that $v(x)\geq 1$ for all $x\in X$.  Therefore $u_n\uparrow \onef$ in $\contX$.  But if $x\in Z(u)$, then $\delta_x(u_n) = u_n(x)=0$ for every $n\in\mathbb{N}$ so that $\delta_x\notin\ordercontc{\contX}$.
\end{rem}

\bibliographystyle{amsplain}
\bibliography{UCompletebib}

\end{document}